\documentclass{amsart}

\usepackage{amsfonts}
\usepackage{cases}
\usepackage{mathrsfs}
\usepackage{bbm}
\usepackage{amssymb}
\usepackage{txfonts}
\usepackage{amscd}
\usepackage{amsfonts,latexsym,amsmath,amsthm,amsxtra,mathdots}
\usepackage[bookmarks,colorlinks]{hyperref}
\usepackage[all,cmtip]{xy}
\RequirePackage{amsmath} \RequirePackage{amssymb}
\usepackage{color}
\usepackage{colordvi}
\usepackage{multicol}
\usepackage{tikz}
\usepackage{hyperref}
\usepackage{graphicx}
\usepackage{amsmath}
\usepackage{amsmath,amscd}
\usetikzlibrary{matrix,arrows}

\def\-{^{-1}}

\newcommand{\delete}[1]{}

    \theoremstyle{plain}

\newtheorem{thm}{Theorem}[section]
\newtheorem{lem}[thm]{Lemma}
\newtheorem{prop}[thm]{Proposition}
\newtheorem{cor}[thm]{Corollary}
\newtheorem{rem}[thm]{Remark}

\newtheorem{add}[thm]{Addendum}

    \numberwithin{equation}{section}

\def\Proof{\noindent{\bf Proof}\quad}
\def\qed{\hfill$\square$\smallskip}

\begin{document}

\title{Isomorphism Conjecture for Baumslag-Solitar Groups}

\author{F. Thomas Farrell}
\address{Department of Mathematics, State University of New York at Binghamton, NY, 13902, U.S.A.}
\email{farrell@math.binghamton.edu}
\thanks{The first author was supported in part by NSF Grant DMS 1206622.}

\author{Xiaolei Wu}
\address{Department of Mathematics, State University of New York at Binghamton, NY, 13902, U.S.A.}
\email{xwu@math.binghamton.edu}


\date{September 22, 2013 and, in revised form, April 1, 2014.}


\begin{abstract}
In this paper, we
prove the K- and L-theoretical Isomorphism Conjecture for Baumslag-Solitar groups with coefficients in an additive category.
\end{abstract}

\maketitle

\section*{Introduction}

In this paper, we confirm that the Isomorphism conjecture is true for all Baumslag-Solitar groups. Our main theorem is as follows
\begin{thm}
The  K- and L-theoretical Isomorphism Conjecture is true for every Baumslag-Solitar group  with coefficients in any additive category.
\end{thm}
\begin{rem}
Independently, G. Gandini, S. Meinert and H. R\"uping proved the Isomorphism Conjecture with coefficients in an additive category for the fundamental group of any graph of abelian groups in \cite{GMR}, which includes all Baumslag-Solitar groups.
\end{rem}

 Recall that the Baumslag-Solitar group $BS(m,n)$ is defined by $\langle a,b~|~b^{-1}a^mb = a^n \rangle$ and all the solvable ones are isomorphic to $BS(1,n)$. Note that $BS(m,n) \cong BS(n,m) \cong BS(-m,-n)$. When $m = \pm n$, $BS(m,\pm m)$ are CAT(0) groups, hence one can use \cite{BL} and \cite{W2} to conclude  $BS(m,\pm m)$ satisfies the Isomorphism Conjecture with coefficients in an additive category. Since the construction of a $K(BS(m,\pm m),1)$ space with universal cover a CAT(0) space is enlightening for part of our proof, we construct it here. We start with two cylinders $S^1_1 \times [0,1]$, $S^1_2 \times [0,1]$ with $S^1_1$ is the standard circle with radius $m$, $S^1_2$ has radius $1$. We will glue $S^1_1 \times \{0\}$ to $S^1_2 \times \{0\}$ by a $m$-fold (orientation preserving) local isometry covering map. And we do the same to $S^1_1 \times \{1\}$ and $S^1_2 \times \{1\}$. What we get is a $K(BS(m,m),1)$ space with universal cover CAT(0). If we reverse orientation on the second gluing, we will get a $K(BS(m,-m),1)$  space with universal cover CAT(0).

In this paper, we will actually prove the Isomorphism conjecture for $BS(m,n)$ with finite wreath products and coefficients in an additive category.  The K- and L-theoretical Isomorphism Conjecture with finite wreath products and coefficients in an additive category is known for all solvable Baumslag-Solitar groups, see \cite{FW} and \cite{W2}. Hence we only need to prove the case when $m > |n| >1$. We will abbreviate the K- and L-theoretical  Isomorphism Conjecture with finite wreath products and coefficients in an additive category by FJCw. The idea is to analyze the preimages of cyclic subgroups corresponding to a certain linear representation of $BS(m,n)$ by using Bass-Serre theory, and conclude they satisfy FJCw. Our results relies on previous work on FJCw for CAT(0) groups and virtually solvable groups by Bartels and L\"uck in  \cite{BL}, and Wegner in \cite{W1}, \cite{W2}. Note that FJCw implies the Isomorphism Conjecture with coefficients in an additive category.

\section{Inheritance Properties and Results on FJCw}

We list some inheritance properties and results on FJCw that we will need. For more information about FJCw we refer to \cite{W2} section 2.3.
\begin{prop}\label{qut}
(1) If  a group $G$
satisfies FJCw, then every subgroup $H < G$ satisfies FJCw.\\
(2) If $G_1$ and $G_2$ satisfy FJCw, then their direct product $G_1 \times G_2$ and free product $G_1 \ast G_2$ satisfy FJCw.\\
(3) Let $\{G_i ~|~ i \in I\}$ be a directed system of groups (with not necessarily injective
structure maps). If each $G_i$ satisfies FJCw, then the direct limit $colim_{i \in I} G_i$
satisfies FJCw.\\
(4) Let $\phi : G \rightarrow Q$ be a
group homomorphism. If $Q$ and $\phi^{-1}(C)$ satisfy FJCw for every cyclic subgroup $C < Q$ then G satisfies FJCw.\\
(5) CAT(0) groups satisfy FJCw.\\
(6) Virtually solvable groups satisfy FJCw. In particular ${\mathbb{Z} {[{\frac{1}{mn}} ]}} \rtimes_{\frac{m}{n}} \mathbb{Z}$ satisfies FJCw, where $m>|n| >0$.\\
\end{prop}

Proof of (1) - (4) can be found for example in \cite{W2} section 2.3. (5) is the main result of \cite{BL} and \cite{W1}. (6) is proved in \cite{W2}.

\section{Proof of the main theorem}

In this section, we prove our main theorem. Recall $BS(m,n) = \langle a,b|b^{-1}a^m b = a^n\rangle$. We assume $m > |n| > 1$.

Let $\Gamma(m,n) = {\mathbb{Z} {[{\frac{1}{mn}} ]}} \rtimes_{\frac{m}{n}} \mathbb{Z}$. $\Gamma(m,n)$ is a solvable linear group. There is a map $ \phi: BS(m,n) \rightarrow \Gamma(m,n)$ mapping $a$ to $(1,0)$ and $b$ to $(0,1)$. Note that FJCw is known for $\Gamma(m,n)$ by Proposition \ref{qut} (6). Hence by Proposition \ref{qut} (4) in order to prove FJCw for $BS(m,n)$, we only need to prove FJCw holds for any subgroup $\phi^{-1}(C)$ of $\Gamma(m,n)$  where $C$ is a cyclic subgroup of $\Gamma(m,n)$.

We proceed to analyze $\phi^{-1}(C)$ using Bass-Serre theory. Part of the idea here is  from the proof of Lemma 4.3 in \cite{FL}.

\begin{figure}[h]

		\centering
		\includegraphics[width=1.0\textwidth]{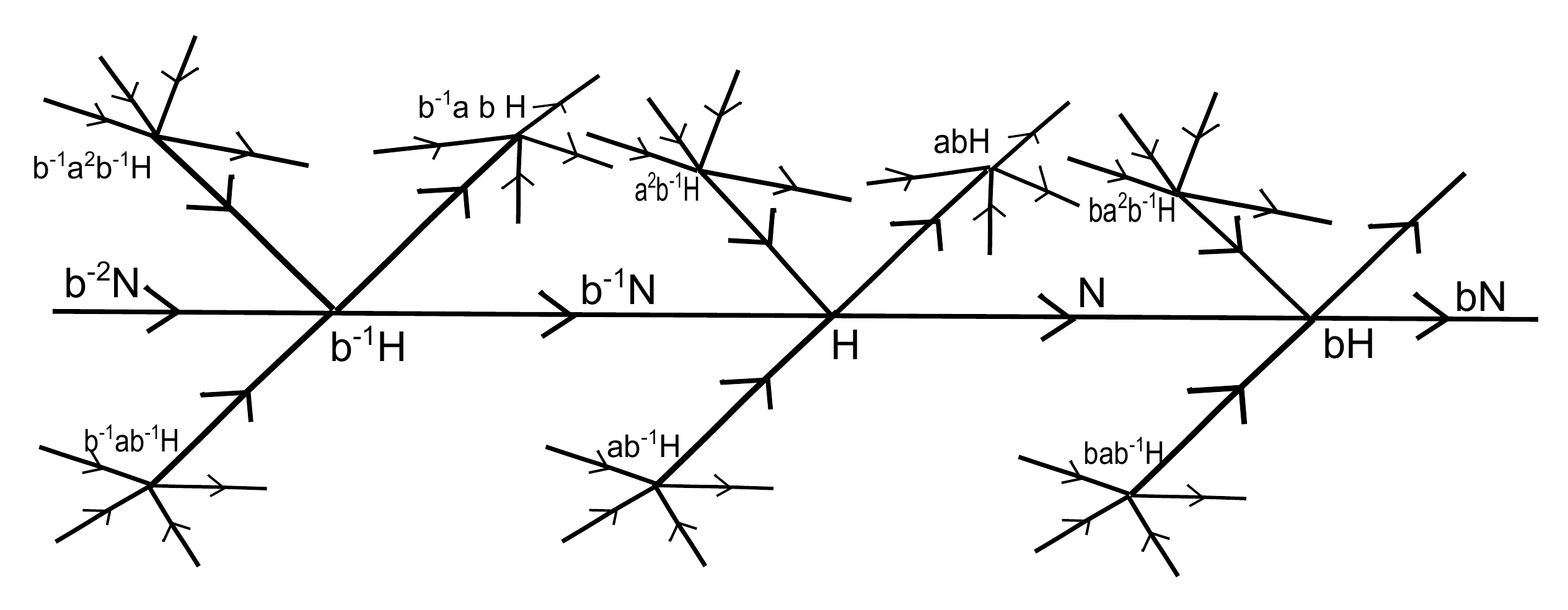}
		\centering
        \caption{T(2,3)}\label{tree2-3}
\end{figure}

Viewing $BS(m,n)$ as an HNN extension, we have an associated  oriented tree $T(m,n)$ and $BS(m,n)$ acts on it without inversion, see for example \cite{DD}, I.3.4. Let $H$ be the cyclic subgroup in $BS(m,n)$ generated by $a$, $N$ be the subgroup generated by $a^m$. Then the vertices of $T(m,n)$ are left cosets $gH$, and edges are the left cosets $gN$ where $g \in BS(m,n)$. The edge $gN$ connects from the tail vertex $gH$ to the head vertex $gbH$. The stabilizer of the vertex $gH$ is the subgroup $gHg^{-1}$, and the stabilizer of the edge $gN$ is $gNg^{-1}$. At each vertex, there are $m$ edges going out and $n$ edges going in. See Figure \ref{tree2-3} for a picture of $T(2,3)$.

\begin{prop} \label{case1}
If $C$ is generated by $(x,y)$ with $y \neq 0$, then $\phi^{-1}(C)$  is a free group. Hence $\phi^{-1}(C)$ satisfies the FJCw.
\end{prop}

\Proof We will show $\phi^{-1}(C)$  is a free group by showing it acts freely on the tree $T(m,n)$ (see for example \cite{DD}, I.4.1). For every vertex $gH$ in $T(m,n)$, its stabilizer under the action of $BS(m,n)$ is the cyclic subgroup $gHg^{-1}$.  Let $\phi(g) = (u,v)$, then for any $k \in \mathbb{Z}$,
$$\phi(ga^kg^{-1}) = \phi(g)  \phi (a^k)\phi(g^{-1}) = (u,v) (k,0)(-{({\frac{n}{m})}^v}u,-v) = ({(\frac{m}{n})}^{v} k ,0) \hspace{1cm} (1) $$
Note that $C$ is generated by $(x,y)$ with $y \neq 0$ and the second coordinate of $(x,y)^n$ is $ny$. Hence $\phi(gHg^{-1}) \cap  C = \{(0,0)\}$. Note that $\phi$ restricted to $H$ is injective, hence it is injective when restricted to $gHg^{-1}$. Therefore $gHg^{-1} \cap \phi^{-1}(C) = \{1\}$. We conclude that for every vertex $gH$, the stabilizer for the action of $\phi^{-1}(C)$  on $T(m,n)$ is trivial. Therefore $\phi^{-1}(C)$ acts freely on $T(m,n)$. Now by Proposition \ref{qut} (5), $\phi^{-1}(C)$ satisfies the FJCw. \qed

\begin{prop} \label{case2}
Let $B = {\mathbb{Z} {[{\frac{1}{mn}} ]}} \rtimes_{\frac{m}{n}} \{0\}$, then $\phi^{-1}(B)$  is  a direct limit of CAT(0) groups. Hence it satisfies FJCw.
\end{prop}

\begin{cor} \label{co}
If $C$ is generated by $(x,0)$, then $\phi^{-1}(C)$ is a subgroup of $\phi^{-1}(B)$. By Proposition \ref{qut} (1), $\phi^{-1}(C)$ also satisfies FJCw.
\end{cor}

We now start the proof of Proposition \ref{case2}. By equation (1) in  the proof of Proposition \ref{case1}, one sees that $\phi(gHg^{-1}) \subset B$ for any $g$. Hence $\phi(gHg^{-1}) \cap B = \phi(gHg^{-1})$, and $gHg^{-1} \cap \phi^{-1}(B) = gHg^{-1} $. Hence vertex or edge stabilizers for the action of $\phi^{-1}(B)$ on $T(m,n)$ are the same as for the $BS(m,n)$ action. In particular, all the vertex and edge stabilizers are isomorphic to $\mathbb{Z}$. Moreover, for every edge $gN$, the stabilizer $G_{gN} \cong \mathbb{Z}$ embeds into the tail vertex stabilizer $G_{gH} \cong \mathbb{Z}$ by multiplying $m$. Correspondingly, $G_{gN} \cong \mathbb{Z}$ embeds into the head vertex stabilizer $G_{gbH} \cong \mathbb{Z}$ by multiplying $n$.

Note if we add one relation $a = 1$ to the group $BS(m,n)$ then the new group is isomorphic to $\mathbb{Z}$. Hence there is a surjective homomorphism $f: BS(m,n) \rightarrow \mathbb{Z}$. If we denote the projection from ${\mathbb{Z} {[{\frac{1}{mn}} ]}} \rtimes_{\frac{m}{n}} \mathbb{Z}$ to $\mathbb{Z}$ by $p$, then $f = p\circ \phi$.
And if $\beta \in \phi^{-1}(B)$, then $f(\beta) = 0$. The action of $\beta$ on $T(m,n)$ has the following key property.

\begin{lem}\label{loo}
For any $\beta \in \phi^{-1}(B)$, the oriented geodesic connecting $gH$ and $\beta gH$ contains an even number of edges, half of them are compatible with the orientation on $T(m,n)$, while the other half are oppositely oriented.
\end{lem}

\Proof We first define an ``algebraic distance" function $D: T(m,n) \times T(m,n) \rightarrow \mathbb{R}$. Note the tree $T(m,n)$ has a standard metric with edge length $1$. For any two points $P,Q \in T(m,n)$, there is a unique oriented geodesic connecting $P$ to $Q$. $D(P,Q)$ is the distance from $P$ to $Q$ counted with signs. In more detail, when the geodesic coincide with tree orientation, it contributes positively; otherwise negatively. For example, in Figure \ref{tree2-3}, $D(H, b^nH) = n$, $D(H, ba^2b^{-1}H) =0$. Note $D(P,Q) = -D(Q,P)$, $D(P,Q) + D(Q,R) = D(P,R)$ for any $P,Q,R \in T(m,n)$.

Since $BS(m,n)$ acts on $T(m,n)$ via orientation preserving isometries, $D(gP,gQ) = D(P,Q)$, for any $g \in BS(m,n)$. In particular, $D(H,gH) = D(aH,agH) = D(H,agH)$ since $a$ is in the stabilizer of the vertex  $H$. And $D(H, gH) = D(bH,bgH)$, hence $D(H, bgH) = D(H,bH) + D(bH,bgH) = 1+ D(H,gH)$.
Note for any $\beta \in \phi^{-1}(B)$, $f(\beta) = 0$. Without loss of generality, we can write $\beta$ as a word $a^{s_1} b^{t_1} a^{s_2}b^{t_2} \dots a^{s_h}b^{t_h}$, $s_i,t_j \in \mathbb{Z}$. Then $t_1+t_2 + \cdots t_h = 0$. Therefore,
$$D(H, \beta H) = D(H, a^{s_1} b^{t_1} a^{s_2}b^{t_2} \cdots a^{s_h}b^{t_h}H)$$
$\hspace*{60mm} = D(H,b^{t_1} a^{s_2}b^{t_2} \cdots a^{s_h}b^{t_h}H)$\\
$\hspace*{60mm} = t_1 + D(H,a^{s_2}b^{t_2} \cdots a^{s_h}b^{t_h}H)$\\
$\hspace*{60mm}\cdots$\\
$\hspace*{60mm} = t_1+t_2 + \cdots t_h$\\
$\hspace*{60mm} =0$ \\
In general, $D(gH, \beta gH) = D(H,g^{-1} \beta g H) = 0$ since $g^{-1} \beta g \in \phi^{-1}(B)$. The lemma now follows. \qed

\begin{add}\label{add}
Lemma \ref{loo} is true for any oriented path $\alpha$ connecting $gH$ and $\beta gH$.
\end{add}

To see this just observe that $\alpha$ can be reduced to the oriented geodesic connecting $gH$ and $\beta gH$ in a finite number of simple moves of the form
$$\alpha = \alpha_1 e (-e) \alpha_2 \longmapsto \alpha_1 \alpha_2$$
where $e$ is an oriented edge and $-e$ is the same edge with the opposite orientation. (And $\alpha_1\alpha_2$, ets. means concatenation of paths.)\\

Summarize what we have so far, the quotient $T(m,n)/ \phi^{-1}(B)$ is an oriented graph of groups with the following properties: every vertex or edge stabilizer is isomorphic to $\mathbb{Z}$, the edge stabilizer embeds into its tail vertex stabilizer by multiplying $m$, and to its head vertex stabilizer by multiplying $n$; moreover, every oriented loop $\alpha$ in $T(m,n)/ \phi^{-1}(B)$ has an even number of edges, where exactly half of them coincide with the orientation. This is a consequence of Addendum \ref{add} and the observation that $\alpha$ lifts to an oriented path $\bar{\alpha}$ in $T(m,n)$ with initial point and end point identified by the action of some element in $\phi^{-1}(B)$.

We will prove that any compact connected subgraph of groups of $T(m,n)/ \phi^{-1}(B)$ is CAT(0), hence $\phi^{-1}(B)$ is a direct limit of CAT(0) groups which satisfies the FJCw by Proposition \ref{qut} (3), (5).

Now let $Y$ be a compact connected subgraph of groups of $T(m,n)/ \phi^{-1}(B)$, and $\Gamma$ be its corresponding group. We will construct a $K(\Gamma,1)$ space with universal cover CAT(0).

Choosing a base vertex $v_0$ in $Y$, since its stabilizer is $\mathbb{Z}$, we construct a cylinder $S^1 \times [0,1]$ with radius $1$. For every edge and every other vertex, we will construct a cylinder $S^1 \times [0,1]$ with radius determined inductively by the following criteria:\\

\textit{\small{The radius of the cylinder corresponding to an edge is $m$ times  the radius of its tail vertex cylinder, and $n$ times the radius of its head vertex cylinder.}}\\

We now explicitly define these radius so that the above criteria are satisfied. We start by defining the radius $r(v)$ of the cylinder corresponding to a vertex $v$ of $Y$ as follows. Let $\alpha$ be an oriented path connecting $v_0$ to $v$. Define the integer $\#(v)$ to be the number of positive (agreeing) edges in $\alpha$ minus the number of negative (disagreeing) edges. To see that $\#(v)$ is well defined, let $\gamma$ be a second such path. Then the concatenation $\alpha(-\gamma)$ is an oriented loop in  $T(m,n)/ \phi^{-1}(B)$ which (as observed above) must have the same number of $+$ edges as $-$ edges. Now let
$$r(v) = (\frac{m}{n})^{(\#(v))}$$
 and then the radius $r(e)$ of an edge $e$ is defined by
 $$r(e) = n (\frac{m}{n})^{(\#(v))}$$
 where $v$ is the head vertex of $e$.


Now we glue all these cylinders together by local isometries according to the graph $Y$. We glue them together using the following method.\\

\textit{\small{ Let $S^1 \times [0, 1]$ be a cylinder corresponding to an edge, denote its tail vertex cylinder by $S_t^1 \times [0,1]$ and its head vertex cylinder by $S_h^1 \times [0,1]$. We glue $S^1 \times \{0\}$ to $S_t^1 \times \{1\}$ by a $m$-fold  locally isometric covering projection (orientation preserving if $m>0$ orientation reversing if $m<0$). And we glue $S^1 \times \{1\}$ to $S_h^1 \times \{0\}$ by a $n$-fold  locally isometric covering projection (orientation preserving if $n>0$ orientation reversing if $n<0$).}}\\

Therefore we have constructed a 2-complex which is a model for $K(\Gamma, 1)$. It is now easy to check that its universal cover is a CAT(0) space; cf. \cite{Da}, p502, Theorem I.2.7. This completes the proof of Proposition \ref{case2}. \qed

Combining Proposition \ref{case1} and Corollary \ref{co}, we prove FJCw for $BS(m,n)$ by using Proposition \ref{qut} (4).

\bibliographystyle{amsplain}

\end{document}